\documentclass{ws-procs9x6}

\newcommand{\R}{\mathrm{I\!R}}
\newcommand{\N}{\mathrm{I\!N}}
\newcommand{\HH}{\mathrm{I\!H}}

\newcommand{\K}{\mathrm{I\!K}}
\newcommand{\PP}{\mathrm{I\!P}}

\newcommand{\Z}{\mathchoice {\hbox{$\sf\textstyle Z\kern-0.4em
Z$}}{\hbox{$\sf\textstyle Z\kern-0.4em Z$}}{\hbox{$\sf\scriptstyle
Z\kern-0.3em Z$}}{\hbox{$\sf\scriptscriptstyle Z\kern-0.2em Z$}}}

\newcommand{\Q}{\mathchoice {\setbox0=\hbox{$\displaystyle\rm
Q$}\hbox{\raise0.15\ht0\hbox to0pt{\kern0.4\wd0\vrule
height0.8\ht0\hss}\box0}}{\setbox0=\hbox{$\textstyle\rm
Q$}\hbox{\raise0.15\ht0\hbox to0pt{\kern0.4\wd0\vrule
height0.8\ht0\hss}\box0}}{\setbox0=\hbox{$\scriptstyle\rm
Q$}\hbox{\raise0.15\ht0\hbox to0pt{\kern0.4\wd0\vrule
height0.7\ht0\hss}\box0}}{\setbox0=\hbox{$\scriptscriptstyle\rm
Q$}\hbox{\raise0.15\ht0\hbox to0pt{\kern0.4\wd0\vrule
height0.7\ht0\hss}\box0}}}

\newcommand{\OO}{\mathchoice {\setbox0=\hbox{$\displaystyle\rm
O$}\hbox{\hbox to0pt{\kern0.4\wd0\vrule
height0.9\ht0\hss}\box0}}{\setbox0=\hbox{$\textstyle\rm O$}\hbox{\hbox
to0pt{\kern0.4\wd0\vrule
height0.9\ht0\hss}\box0}}{\setbox0=\hbox{$\scriptstyle\rm O$}\hbox{\hbox
to0pt{\kern0.4\wd0\vrule
height0.9\ht0\hss}\box0}}{\setbox0=\hbox{$\scriptscriptstyle\rm
O$}\hbox{\hbox to0pt{\kern0.4\wd0\vrule height0.9\ht0\hss}\box0}}}

\newcommand{\id}{\mathrm{id}}

\newcommand{\SO}{\mathrm{SO}}
\newcommand{\End}{\mathrm{End}}

\newcommand{\rk}{\mathrm{rk}}

\newcommand{\vkap}{\varkappa}

\newcommand{\Menge}[2]{\{\,#1\,|\,#2\,\}}

\newcommand {\g}[2]{\langle #1,#2\rangle}

\newcommand {\operp}{\mathbin{\mbox{$\ominus\raisebox{2.9pt}
 {\hskip-0.42em\hbox{\vrule height0.7ex width0.02em}\hskip0.42em }$}}}

\newcommand{\Ug}{\mathrm{U}}
\newcommand{\SU}{\mathrm{SU}}
\newcommand{\Sp}{\mathrm{Sp}}

\newcommand{\Spin}{\mathrm{Spin}}

\newcommand{\Fix}{\mathop{\mathrm{Fix}}\nolimits}

\newcommand{\RP}{\ensuremath{\R\mathrm{P}}}
\newcommand{\CP}{\ensuremath{\C\mathrm{P}}}
\newcommand{\HP}{\ensuremath{\HH\mathrm{P}}}
\newcommand{\OP}{\ensuremath{\OO\mathrm{P}}}

\newcommand{\bbS}{\mathbb{S}}

\newcommand{\Sph}{\bbS}

\newcommand{\C}{\mathchoice {\setbox0=\hbox{$\displaystyle\rm
C$}\hbox{\hbox to0pt{\kern0.4\wd0\vrule
height0.95\ht0\hss}\box0}}{\setbox0=\hbox{$\textstyle\rm C$}\hbox{\hbox
to0pt{\kern0.4\wd0\vrule
height0.95\ht0\hss}\box0}}{\setbox0=\hbox{$\scriptstyle\rm C$}\hbox{\hbox
to0pt{\kern0.4\wd0\vrule
height0.95\ht0\hss}\box0}}{\setbox0=\hbox{$\scriptscriptstyle\rm
C$}\hbox{\hbox to0pt{\kern0.4\wd0\vrule height0.95\ht0\hss}\box0}}}

\begin{document}

\title{Totally geodesic submanifolds in Riemannian symmetric spaces}

\author{Sebastian Klein}

\address{Universit\"at Mannheim, Lehrstuhl f\"ur Mathematik III, \\
Seminargeb\"aude A5, 68131 Mannheim, Germany \\
E-Mail: mail@sebastian-klein.de}

\begin{abstract}
In the first part of this expository article, the most important constructions and classification results concerning totally geodesic submanifolds
in Riemannian symmetric spaces are summarized. In the second part, I describe the results of my classification of the totally geodesic submanifolds
in the Riemannian symmetric spaces of rank 2.

To appear in the Proceedings volume for the conference \emph{VIII International Conference on Differential Geometry}, which took place in Santiago de Compostela
in July 2008.
\end{abstract}

\bodymatter

\section{Totally geodesic submanifolds}
\label{Se:intro}

A submanifold \,$M'$\, of a Riemannian manifold \,$M$\, is called \emph{totally geodesic}, if every geodesic of \,$M'$\, 
is also a geodesic of \,$M$\,. In this article, we will discuss totally geodesic submanifolds in Riemannian symmetric spaces;
in such spaces, a connected, complete submanifold is totally geodesic if and only if it is a symmetric subspace.

There are several important construction principles for totally geodesic submanifolds in Riemannian symmetric spaces \,$M$\,. First, we note that the connected
components of the fixed point set of any isometry \,$f$\, of \,$M$\, are totally geodesic submanifolds (this is in fact true in any Riemannian manifold,
see Ref.~\refcite{Kobayashi:1972}, Theorem~II.5.1, p.~59). This construction principle is especially important in the case where \,$f$\, is involutive
(i.e.~\,$f\circ f=\id_M$\,); the totally geodesic submanifolds resulting in this case are called \emph{reflective} submanifolds;
they have been studied extensively, for example by \textsc{Leung} (see below).

Further constructions of totally geodesic submanifolds in Riemannian symmetric spaces of compact type \,$M$\, were introduced by \textsc{Chen} and \textsc{Nagano}%
\cite{Chen/Nagano:totges2-1978, Chen:1987}:
For \,$p\in M$\,, the connected components \,$\neq \{p\}$\, of the fixed point set of the geodesic reflection
of \,$M$\, at \,$p$\, are called \emph{polars} or \emph{$M_+$-submanifolds} of \,$M$\,; note that they are in particular reflective submanifolds of \,$M$\,.
A \emph{pole} of \,$M$\, is a polar which is a singleton. It has been shown by Chen/Nagano\cite{Chen/Nagano:totges2-1978} that for every polar \,$M_+$\, of \,$M$\,
and every \,$q \in M_+$\, there exists another reflective submanifold \,$M_-$\, of \,$M$\, with \,$q \in M_-$\, and \,$T_qM_- = (T_qM_+)^\perp$\,; \,$M_-$\, is
called a \emph{meridian} or \emph{$M_-$-submanifold} of \,$M$\,. For \,$p_1,p_2 \in M$\,, a point \,$q\in M$\, is called a \emph{center point} between 
\,$p_1$\, and \,$p_2$\, if there exists a geodesic joining \,$p_1$\, with \,$p_2$\, so that \,$q$\, is the middle point on that geodesic. If \,$p_2$\, 
is a pole of \,$p_1$\,, then the set \,$C(p_1,p_2)$\, of center points between \,$p_1$\, and \,$p_2$\, is called the \emph{centrosome} of \,$p_1$\, and \,$p_2$\,;
its connected components are totally geodesic submanifolds of \,$M$\, (see Ref.~\refcite{Chen:1987}, Proposition~5.1). 

Moreover, every symmetric space of compact type can be embedded in its transvection group as a totally geodesic submanifold: Let \,$M=G/K$\, be such a space,
then there exists an involutive automorphism \,$\sigma$\, of \,$G$\, so that \,$\Fix(\sigma)^0 \subset K \subset \Fix(\sigma)$\,. Because of this property,
the \emph{Cartan map}
$$ f: G/K \to G, \; gK \mapsto \sigma(g)\cdot g^{-1} $$
is a well-defined covering map onto its image, which turns out to be a totally geodesic submanifold of \,$G$\,. If \,$M$\, is a ``bottom space'', i.e.~there
exists no non-trivial symmetric covering map with total space \,$M$\,, then we have \,$K = \Fix(\sigma)$\,, and therefore \,$f$\, is an embedding. In this
setting \,$f$\, is called the \emph{Cartan embedding} of \,$M$\,. 

\medskip

It is a significant and interesting problem to determine all totally geodesic submanifolds in a given symmetric space. Because totally 
geodesic submanifolds are rigid (i.e.~if \,$M'_1,M'_2$\, are connected, complete totally geodesic submanifolds of \,$M$\,
with \,$p\in M'_1\cap M'_2$\, and \,$T_pM_1' = T_pM_2'$\,, then we already have \,$M_1=M_2$\,), they can be classified
by determining those linear subspaces \,$U \subset T_pM$\, which occur as tangent spaces of totally geodesic submanifolds of \,$M$\,. 

The elementary answer to the latter problem is the following: 
There exists a totally geodesic submanifold of \,$M$\,
with a given tangent space \,$U\subset T_pM$\, if and only if \,$U$\, is \emph{curvature invariant} (i.e.~we have \,$R(u,v)w \in U$\, for all \,$u,v,w \in U$\,,
denoting by \,$R$\, the Riemannian curvature tensor of \,$M$\,). 

Therefore the classification of totally geodesic submanifolds of \,$M$\, reduces to the purely algebraic problem of the
classification of curvature invariant subspaces of \,$T_pM$\,. However, because of the algebraic complexity of the curvature tensor, classifying
the curvature invariant subspaces is by no means an easy task, and therefore also the classification of totally geodesic submanifolds
remains a signification problem.

This problem has been solved
for the Riemannian symmetric spaces of rank~1 by \textsc{Wolf} in Ref.~\refcite{Wolf:1963-elliptic}, \S 3. 
\textsc{Chen}/\textsc{Nagano} claimed a classification for the complex quadrics (which are symmetric spaces
of rank 2) in Ref.~\refcite{Chen/Nagano:totges1-1977}, and then for all symmetric spaces of rank 2 
in Ref.~\refcite{Chen/Nagano:totges2-1978}, using their construction of polars and meridians described above.
However, it turns out that their classifications are incorrect: 
For several spaces of rank 2, totally geodesic submanifolds have been missed, and also some other details are faulty. In my papers 
Refs.~\refcite{Klein:2007-claQ, Klein:2007-tgG2, Klein:2007-Satake, Klein:2008-tgex} I discuss these shortcomings and give a full classification
of the totally geodesic submanifolds in all irreducible symmetric spaces of rank 2; Section~\ref{Se:rk2} of the present exposition
contains a summary of these results. For symmetric spaces of rank \,$\geq 3$\,, the full classification problem is still open.

However, there are several results concerning the classification of special classes of totally geodesic submanifolds. Probably the most significant
result of this kind is the classification of reflective submanifolds in all Riemannian symmetric spaces due to \textsc{Leung};
his results are found in final form in Ref.~\refcite{Leung:reflective-1979}, but also see Refs.~\refcite{Leung:reflective-1974, Leung:reflective-errata-1975}. 
Another important problem of this kind is the classification of the totally geodesic submanifolds \,$M'$\,
of \,$M=G/K$\, with maximal rank (i.e.~\,$\rk(M') = \rk(M)$\,); this problem has been solved for the symmetric spaces with \,$\rk(M)=\rk(G)$\,
by \textsc{Ikawa}/\textsc{Tasaki}\cite{Ikawa/Tasaki:2000}, and then for all irreducible symmetric spaces by \textsc{Zhu}/\textsc{Liang}\cite{Zhu/Liang:2004}. 

Further important classification results concern Hermitian symmetric spaces \,$M$\,: In them, the complex totally geodesic submanifolds have
been classified by \textsc{Ihara}\cite{Ihara:tgcomplex-1967}. Moreover, the real forms of \,$M$\, (i.e.~the totally real, totally geodesic submanifolds
\,$M'$\, of \,$M$\, with \,$\dim_{\R}(M') = \dim_{\C}(M)$\,) are all reflective; due to this fact \textsc{Leung} was able to derive a classification
of the real forms of all Hermitian symmetric spaces from his classification of reflective submanifolds.\cite{Leung:realforms-1979}.

Finally we mention a result by \textsc{Wolf} concerning totally geodesic submanifolds \,$M'$\, of (real, complex or quaternionic) Grassmann manifolds
\,$G_r(\K^n)$\, with the property that any two distinct elements of \,$M'$\, have zero intersection, regarded as \,$r$-dimensional
subspaces of \,$\K^n$\,. Wolf showed\cite{Wolf:1963-spheres, Wolf:1963-elliptic}
that any such totally geodesic submanifold is isometric either to a sphere or to a projective space over \,$\R$\,,
\,$\C$\, or \,$\HH$\,; he was also able to describe embeddings for these submanifolds explicitly and to calculate their maximal dimension
in dependence of \,$r$\, and \,$n$\,. 

\section{Maximal totally geodesic submanifolds in the Riemannian symmetric spaces of rank 2}
\label{Se:rk2}

In the following, I list the isometry types corresponding to all the congruence classes of 
totally geodesic submanifolds which are maximal (i.e.~not strictly
contained in another connected totally geodesic submanifold) in all Riemannian symmetric spaces of rank 2.
In many cases I also briefly describe totally geodesic embeddings corresponding
to these submanifolds. This is a summary of my work in Refs.~\refcite{Klein:2007-claQ, Klein:2007-tgG2, Klein:2007-Satake, Klein:2008-tgex}, where it is proved
that the lists given here are complete, and where the totally geodesic embeddings are described in more detail.

The invariant Riemannian metric of an irreducible Riemannian symmetric space is unique only up to a positive constant. In the sequel, we use the following notations
to describe the metric which is induced on the totally geodesic submanifolds: 
For \,$\ell \in \N$\, and \,$r > 0$\, we denote by \,$\Sph^\ell_r$\, the \,$\ell$-dimensional sphere of radius \,$r$\,, and for \,$\vkap > 0$\,
we denote by \,$\RP^\ell_\vkap$\,, \,$\CP^\ell_\vkap$\,, \,$\HP^\ell_\vkap$\, and \,$\OP^2_\vkap$\, 
the respective projective spaces, their metric being scaled in such a way that the \emph{minimal}
sectional curvature is \,$\vkap$\,. (\,$\RP^\ell_\vkap$\, is then of constant sectional curvature \,$\vkap$\,, \,$\CP^\ell_\vkap$\, is of constant holomorphic
sectional curvature \,$4\vkap$\,, and we have the inclusions \,$\RP^\ell_\vkap \subset \CP^\ell_\vkap \subset \HP^\ell_\vkap$\, of totally geodesic submanifolds). 
For symmetric spaces of rank 2, we describe the appropriate metric by stating the length \,$a$\, of the shortest restricted root of the space as a subscript
\,${}_{\mathrm{srr}=a}$\,. 
For the three infinite families of Grassmann manifolds \,$G_2^+(\R^n)$\,, \,$G_2(\C^n)$\, and \,$G_2(\HH^n)$\,,
we also use the notation \,${}_{\textrm{srr}=1*}$\, to denote the metric scaled in such a way that the shortest root
occurring \emph{for large \,$n$\,} has length \,$1$\,, disregarding the fact that this root might vanish for certain small values of \,$n$\,.

The spaces in which the totally geodesic submanifolds are classified below are always taken with \,${}_{\mathrm{srr}=1*}$\, (for the Grassmann manifolds)
or \,${}_{\mathrm{srr}=1}$\, (for all others).

\subsection[$G_2^+(\R^{n+2})$]{$\boldsymbol{G_2^+(\R^{n+2})}$}
\label{SSe:rk2:G2Rn}

%See Ref.~\refcite{Klein:2007-claQ} for the classification and the embeddings for most types of totally geodesic submanifolds; Ref.~\refcite{Klein:2007-tgG2}, Section~7 for
%the embedding for the totally geodesic submanifold of type (e).

\begin{enumerate}
\item \,$G_2^+(\R^{n+1})_{\mathrm{srr}=1*}$ \\
The linear isometric embedding \,$\R^{n+1} \to \R^{n+2},\;(x_1,\dotsc,x_{n+1}) \mapsto (x_1,\dotsc,x_{n+1},0)$\, induces a totally geodesic, isometric embedding
\,$G_2^+(\R^{n+1}) \to G_2^+(\R^{n+2})$\,. 

\item \,$\Sph^n_{r=1}$\, \\
Fix a unit vector \,$v_0 \in \R^{n+2}$\,, and let \,$\Sph := \Menge{v \in \R^{n+2}}{\g{v}{v_0}=0,\;\|v\|=1} \cong \Sph^n_{r=1}$\,. Then 
the map \,$\Sph \to G_2^+(\R^{n+2}),\;v \mapsto \R v \oplus \R v_0$\, is a totally geodesic, isometric embedding.

\item \,$(\Sph^\ell_{r=1} \times \Sph^{\ell'}_{r=1})/\Z_2$\,, where \,$\ell+\ell'=n$\, \\
The map
{\footnotesize
\begin{align*}
\Sph^\ell_{r=1} \times \Sph^{\ell'}_{r=1} & \to G_2^+(\R^{n+2}) \\
((x_0,\dotsc,x_\ell),(y_0,\dotsc,y_{\ell'})) & \mapsto \R\,(x_0,\dotsc,x_\ell,0,\dotsc,0) \oplus \R\,(0,\dotsc,0,y_0,\dotsc,y_{\ell'}) 
\end{align*}

}
is a totally geodesic, isometric immersion, and a two-fold covering map onto its image in \,$G_2^+(\R^{n+2})$\,. 

\item For \,$n\geq 4$\, even: \,$\CP^{n/2}_{\vkap=1/2}$\, \\
Let us fix a complex structure \,$J$\, on \,$\R^{n+2}$\,. Then the complex-1-dimensional linear subspaces of \,$(\R^{n+2},J)$\, are in particular
real-2-dimensional oriented linear subspaces of \,$\R^{n+2}$\,. Therefore the complex projective space \,$\PP \cong \CP^{n/2}$\, over \,$(\R^{n+2},J)$\, 
is contained in \,$G_2^+(\R^{n+2})$\,; it turns out to be a totally geodesic submanifold.

\item For \,$n=2$\,: \,$\CP^1_{\vkap=1/2} \times \RP^1_{\vkap=1/2}$\, \\
The image of the Segr\'e embedding \,$\CP^1 \times \CP^1 \to \CP^3$\, (see for example Ref.~\refcite{Shafarevich:1994}, p.~55f.) 
is a 2-dimensional complex quadric in \,$\CP^3$\,; such a quadric is isometric
to \,$G_2^+(\R^4)$\,. Thereby we see that \,$G_2^+(\R^4)$\, is isometric to \,$\CP^1_{\vkap=1/2} \times \CP^1_{\vkap=1/2}$\,. Let \,$C$\, be the trace of
a (closed) geodesic in \,$\CP^1_{\vkap=1/2}$\,; then \,$C$\, is isometric to \,$\RP^1_{\vkap=1/2}$\,, and \,$\CP^1_{\vkap=1/2}\times C$\, is a totally geodesic
submanifold of \,$\CP^1_{\vkap=1/2}\times \CP^1_{\vkap=1/2} \cong G_2^+(\R^4)$\,. 

\item For \,$n=3$\,: \,$\Sph^2_{r=\sqrt{5}}$\, \\
To describe this totally geodesic submanifold, as well as similar totally geodesic submanifolds occurring in \,$G_2(\C^6)$\, and \,$G_2(\HH^7)$\, (see
Sections~\ref{SSe:rk2:G2Cn}(g) and \ref{SSe:rk2:G2Hn}(f) below), we note that there is exactly one irreducible, 14-dimensional, quaternionic
representation of \,$\Sp(3)$\, (see Ref.~\refcite{Broecker-tom-Dieck:1985}, Chapter~VI, Section~(5.3), p.~269ff.). It can be constructed as follows: 
The vector representation of \,$\Sp(3)$\, on \,$\C^6$\, induces a representation of \,$\Sp(3)$\, on \,$\bigwedge^3 \C^6$\,. This
20-dimensional representation decomposes into two irreducible components: One, 6-dimensional, is equivalent to the vector representation of \,$\Sp(3)$\,;
the other, acting on a 14-dimensional linear space \,$V$\,, is the irreducible representation we are interested in.

It turns out that the restriction of the representation of \,$\Sp(3)$\, on \,$V$\, to an \,$\SO(3)$\, embedded in \,$\Sp(3)$\, in the canonical way,
is a real representation, and that in any real form \,$V_{\R}$\, of \,$V$\,, two linear independent vectors are left invariant. By splitting off the subspace
of \,$V'$\, spanned by these vectors, we get a real-5-dimensional representation \,$V_{\R}'$\, of \,$\SO(3)$\,, which turns out to be irreducible
(and equivalent to the Cartan representation \,$\SO(3) \times \End_+(\R^3)_0 \to \End_+(\R^3)_0, \; (B,X) \mapsto BXB^{-1}$\,). It turns out that
the corresponding action of \,$\SO(3)$\, on \,$G_2^+(V_{\R}') \cong G_2^+(\R^5)$\, has exactly one totally geodesic orbit; this orbit is isometric to \,$\Sph^2_{r=\sqrt{5}}$\,. 
\end{enumerate}

\subsection[$G_2(\C^{n+2})$]{$\boldsymbol{G_2(\C^{n+2})}$}
\label{SSe:rk2:G2Cn}

% See Ref.~\refcite{Klein:2007-tgG2}, Section~7.

\begin{enumerate}
\item \,$G_2(\C^{n+1})_{\mathrm{srr}=1*}$\, \\
The linear isometric embedding \,$\C^{n+1} \to \C^{n+2},\;(z_1,\dotsc,z_{n+1}) \mapsto (z_1,\dotsc,z_{n+1},0)$\, induces a totally geodesic, isometric embedding
\,$G_2(\C^{n+1}) \to G_2(\C^{n+2})$\,. 

\item \,$G_2(\R^{n+2})_{\mathrm{srr}=1*}$\, \\
The map \,$G_2(\R^{n+2}) \to G_2(\C^{n+2}),\;\Lambda \mapsto \Lambda\oplus i\Lambda$\, is a totally geodesic, isometric embedding.

\item \,$\CP^n_{\vkap=1}$\, \\
Fix a unit vector \,$v_0 \in \C^{n+2}$\,. Then the map \,$\CP((\C v_0)^\perp) \to G_2(\C^{n+2}),\; \C v \mapsto \C v \oplus \C v_0$\,
is a totally geodesic, isometric embedding.

\item \,$\CP^\ell_{\vkap=1} \times \CP^{\ell'}_{\vkap=1}$\, with \,$\ell+\ell'=n$\, \\
Let \,$\C^{n+2} = W \operp W'$\, be a splitting of \,$\C^{n+2}$\, into complex-linear subspaces of dimension \,$\ell+1$\, resp.~\,$\ell'+1$\,. Then
\,$\CP(W) \times \CP(W') \to G_2(\C^{n+2}),\; (\C v,\C v') \mapsto \C v \oplus\C v'$\, is a totally geodesic, isometric embedding.

\item For \,$n$\, even: \,$\HP^{n/2}_{\vkap=1/2}$\, \\
Let us fix a quaternionic structure \,$\tau$\, on \,$\C^{n+2}$\, (i.e.~\,$\tau: \C^{n+2}\to \C^{n+2}$\, is anti-linear with \,$\tau^2 = -\id$\,).
Then the quaternionic-1-dimensional linear subspaces of \,$(\C^{n+2},\tau)$\, are in particular
complex-2-dimensional linear subspaces of \,$\C^{n+2}$\,. Therefore the quaternionic projective space \,$\PP \cong \HP^{n/2}$\, over \,$(\C^{n+2},\tau)$\, 
is contained in \,$G_2(\C^{n+2})$\,; it turns out that \,$\PP$\, is a totally geodesic submanifold of \,$G_2(\C^{n+2})$\,. 

\item For \,$n=2$\,: \,$G_2^+(\R^5)_{\mathrm{srr}=\sqrt{2}}$\, and \,$(\Sph^3_{r=1/\sqrt{2}} \times \Sph^1_{r=1/\sqrt{2}})/\Z_2$\, \\
Note that \,$G_2(\C^4)_{\mathrm{srr}=1*}$\, is isometric to \,$G_2^+(\R^6)_{\mathrm{srr}=\sqrt{2}}$\,. 
This isometry can be exhibited via the Pl\"ucker map \,$G_2(\C^m) \to \CP(\bigwedge^2 \C^m),\;
\C u \oplus \C v \mapsto \C(u\wedge v)$\,, which is an isometric embedding for any \,$m$\,; for \,$m=4$\, its image in \,$\CP(\bigwedge^2 \C^4)\cong \CP^5$\,
turns out to be a 4-dimensional complex quadric; such a quadric is isomorphic to \,$G_2^+(\R^6)$\,. Thus \,$G_2(\C^4)$\, is isometric to \,$G_2^+(\R^6)$\,,
hence its maximal totally geodesic submanifolds are those given in Section~\ref{SSe:rk2:G2Rn} for \,$n=4$\,, namely: \,$G_2^+(\R^5)_{\mathrm{srr}=\sqrt{2}}$\,, 
\,$(\Sph^3_{r=1/\sqrt{2}} \times \Sph^1_{r=1/\sqrt{2}})/\Z_2$\,, \,$(\Sph^2_{r=1/\sqrt{2}} \times \Sph^2_{r=1/\sqrt{2}})/\Z_2$\,, \,$\Sph^4_{r=1/\sqrt{2}}$\,, \,$\CP^2_{\vkap=1}$\,. The first
two of these submanifolds are those which are listed under this point; the remaining submanifolds have already been listed above (note that \,$(\Sph^2\times \Sph^2)/\Z_2$\, 
and \,$\Sph^4$\, are isometric to \,$G_2(\R^4)$\, and \,$\HP^1$\,, respectively). 

\item For \,$n=4$\,: \,$\CP^2_{\vkap=1/5}$\, \\
Let us consider the 14-dimensional quaternionic, irreducible representation \,$V$\, of \,$\Sp(3)$\, described in Section~\ref{SSe:rk2:G2Rn}(f). The restriction
of that representation to a \,$\SU(3)$\, canonically embedded in \,$\Sp(3)$\, leaves a totally complex 6-dimensional linear subspace \,$V_{\C}$\, of \,$V$\,
invariant; the resulting 6-dimensional representation \,$V_{\C}$\, of \,$\SU(3)$\, is irreducible. It turns out that the induced action of \,$\SU(3)$\,
on \,$G_2(V_{\C}) \cong G_2(\C^6)$\, has exactly one totally geodesic orbit; this orbit is isometric to \,$\CP^2_{\vkap=1/5}$\,. 
\end{enumerate}

\subsection[$G_2(\HH^{n+2})$]{$\boldsymbol{G_2(\HH^{n+2})}$}
\label{SSe:rk2:G2Hn}

% See Ref.~\refcite{Klein:2007-tgG2}: Section~5 for the classification, Section~6 for the embeddings.

\begin{enumerate}
\item \,$G_2(\HH^{n+1})_{\mathrm{srr}=1*}$\, \\
The linear isometric embedding \,$\HH^{n+1} \to \HH^{n+2},\;(q_1,\dotsc,q_{n+1}) \mapsto (q_1,\dotsc,q_{n+1},0)$\, induces a totally geodesic, isometric embedding
\,$G_2(\HH^{n+1}) \to G_2(\HH^{n+2})$\,. 

\item \,$G_2(\C^{n+2})_{\mathrm{srr}=1*}$\, \\
We fix two orthogonal imaginary unit quaternions \,$i$\, and \,$j$\,, and let \,$\C = \R \oplus \R i$\,. Then the map
\,$G_2(\C^{n+2}) \to G_2(\HH^{n+2}), \; \Lambda \mapsto \Lambda \oplus \Lambda j$\, is a totally geodesic, isometric embedding.

\item \,$\HP^n_{\vkap=1}$\, \\
Fix a unit vector \,$v_0 \in \HH^{n+2}$\,. Then the map \,$\HP((v_0\HH)^\perp) \to G_2(\HH^{n+2}),\; v\HH \mapsto v\HH \oplus v_0\HH$\,
is a totally geodesic, isometric embedding.

\item \,$\HP^\ell_{\vkap=1} \times \HP^{\ell'}_{\vkap=1}$\, with \,$\ell+\ell'=n$\, \\
Let \,$\HH^{n+2} = W \operp W'$\, be a splitting of \,$\HH^{n+2}$\, into quaternionic-linear subspaces of dimension \,$\ell+1$\, resp.~\,$\ell'+1$\,. Then
\,$\HP(W) \times \HP(W') \to G_2(\HH^{n+2}),\; (v\HH,v'\HH) \mapsto v\HH\oplus v'\HH$\, is a totally geodesic, isometric embedding.

\item For \,$n=2$\,: \,$\Sp(2)_{\mathrm{srr}=\sqrt{2}}$\, and \,$(\Sph^5_{r=1/\sqrt{2}} \times \Sph^1_{r=1/\sqrt{2}})/\Z_2$\, \\
Let \,$U \in G_2(\HH^4)$\, be given, then \,$U^\perp$\, is the only pole corresponding to \,$U$\, in \,$G_2(\HH^4)$\,. The centrosome between this pair
of poles is a totally geodesic submanifold of \,$G_2(\HH^4)$\, which is isometric to \,$\Sp(2)$\,. This \,$\Sp(2)$\, is also a reflective submanifold of
\,$G_2(\HH^4)$\,, the complementary reflective submanifold is isometric to \,$(\Sph^5_{r=1/\sqrt{2}} \times \Sph^1_{r=1/\sqrt{2}})/\Z_2$\,.

\item For \,$n=5$\,: \,$\HP^2_{\vkap=1/5}$\, \\
We again consider the irreducible, quaternionic 14-dimensional representation \,$V$\, of \,$\Sp(3)$\, introduced in Section~\ref{SSe:rk2:G2Rn}(f);
we now view \,$V$\, as a quaternionic-7-dimensional linear space.
The representation of \,$\Sp(3)$\, on \,$V$\, induces an action of \,$\Sp(3)$\, on the
quaternionic 2-Grassmannian \,$G_2(V) \cong G_2(\HH^7)$\,; again it turns out that this action has exactly one totally geodesic orbit, which is isometric
to \,$\HP^2_{\vkap=1/5}$\,. 

\item For \,$n=4$\,: \,$\Sph^3_{r=2\sqrt{5}}$\, \\
According to the present list, two of the maximal totally geodesic submanifolds of the 2-Grassmannian \,$G_2(\HH^7)$\, are isometric to
\,$G_2(\HH^6)$\, and \,$\HP^2_{\vkap=1/5}$\,, respectively. The intersection of these two totally geodesic submanifolds is a totally geodesic
submanifold of \,$G_2(\HH^6)$\,, which turns out to be isometric to \,$\Sph^3_{r=2\sqrt{5}}$\,.
\end{enumerate}

\subsection[$\SU(3)/\SO(3)$]{$\boldsymbol{\SU(3)/\SO(3)}$}
\label{SSe:rk2:AI}

% See Ref.~\refcite{Klein:2007-Satake}, Section~6.

\begin{enumerate}
\item \,$(\Sph^2_{r=1} \times \Sph^1_{r=\sqrt{3}})/\Z_2$\,
%We consider the map 
%$$ \Phi: \Sph^1_{r=\sqrt{12}} \times \SU(2) \to \SU(3),\; (\lambda,B) \mapsto \begin{pmatrix} \tfrac{\lambda}{\sqrt{12}}\,B & 0 \\ 0 & \tfrac{12}{\lambda^{2}} \end{pmatrix} \; , $$
%where we regard \,$S^1_{r=\sqrt{12}}$\, as the circle \,$\Mengegr{z\in\C}{|z|^2=12}$\, in \,$\C$\,. We have \,$\Phi^{-1}(\SO(3)) = (\{\pm\sqrt{12}\} \times \SO(2))
%\,\cup\, (\{\pm \sqrt{12}i\} \times J\cdot \SO(2)) =: K$\, with \,$J := \left( \begin{smallmatrix} i & 0 \\ 0 & -i \end{smallmatrix} \right) \in \SU(2)$\,.
%Hence \,$\Phi$\, induces a map
%$$ \underline{\Phi}: (\Sph^1_{r=\sqrt{12}} \times \SU(2))/K \to \SU(3)/\SO(3) \; , $$
%which turns out to be a totally geodesic, isometric embedding. Moreover, \,$(\Sph^1_{r=\sqrt{12}} \times \SU(2))/K$\, is isometric to 
%\,$(\Sph^1_{r=\sqrt{3}} \times \Sph^2_{r=1})/\Z_2$\,. 

\item \,$\RP^2_{\vkap=1/4}$\, \\
\end{enumerate}

\subsection[$\SU(6)/\Sp(3)$]{$\boldsymbol{\SU(6)/\Sp(3)}$}
\label{SSe:rk2:AII}

% See Ref.~\refcite{Klein:2008-tgex}, Section~4.4.

\begin{enumerate}
\item \,$\HP^2_{\vkap=1/4}$\, 
%This is the polar in \,$\SU(6)/\Sp(3)$\,. 
\item \,$\CP^3_{\vkap=1/4}$\, 
\item \,$\SU(3)_{\mathrm{srr}=1}$\, \\
The map \,$\SU(3) \to \SU(6)/\Sp(3),\;B \mapsto \left( \begin{smallmatrix} B & 0 \\ 0 & B^{-1} \end{smallmatrix} \right) \cdot \Sp(3)$\,
is a totally geodesic embedding of this type.
\item \,$(\Sph^5_{r=1} \times \Sph^1_{r=\sqrt{3}})/\Z_2$\, 
%This is the meridian in \,$\SU(6)/\Sp(3)$\,. 
\end{enumerate}

\subsection[$\SO(10)/\Ug(5)$]{$\boldsymbol{\SO(10)/\Ug(5)}$}
\label{SSe:rk2:DIII}

% See Ref.~\refcite{Klein:2008-tgex}, Section~3.5.

In the descriptions of the embeddings for this symmetric space, we consider both \,$\Ug(5)$\, and \,$\SO(10)$\, as acting on \,$\C^5 \cong \R^{10}$\,; in the latter
case, this action is only \,$\R$-linear.

\begin{enumerate}
\item \,$\CP^4_{\vkap=1}$\, 
%This is a polar in \,$\SO(10)/\Ug(5)$\,. 
\item \,$G_2(\C^5)_{\mathrm{srr}=1}$\, 
%This is a polar in \,$\SO(10)/\Ug(5)$\,. 
\item \,$\CP^3_{\vkap=1} \times \CP^1_{\vkap=1}$\, \\
%This is the meridian in \,$\SO(10)/\Ug(5)$\, corresponding to \,$G_2(\C^5)$\,. It can also be obtained in the following way: 
Let \,$G := \SO(6)\times \SO(4)$\, be canonically embedded in \,$\SO(10)$\, in such a way that its intersection with \,$\Ug(5)$\, is maximal. Then
\,$G / (G\cap \Ug(5))$\, is a totally geodesic submanifold of \,$\SO(10)/\Ug(5)$\, which is isometric to \,$(\SO(6)/\Ug(3)) \,\times\, (\SO(4)/\Ug(2))
\cong \CP^3 \times \CP^1$\,. 
\item \,$G_2^+(\R^8)_{\mathrm{srr}=\sqrt{2}}$\, \\
%This is the meridian in \,$\SO(10)/\Ug(5)$\, corresponding to \,$\CP^4_{\vkap=1}$\,. It can also be obtained in the following way: 
Let \,$G := \SO(8)$\, be canonically embedded in \,$\SO(10)$\, in such a way that its intersection with \,$\Ug(5)$\, is maximal. 
Then \,$G/(G\cap \Ug(5))$\, is a totally geodesic submanifold of \,$\SO(10)/\Ug(5)$\, which is isometric to \,$\SO(8)/\Ug(4) \cong G_2^+(\R^8)$\,. 
\item \,$\SO(5)_{\mathrm{srr}=1}$\, \\
The map \,$\SO(5) \to \SO(10)/\Ug(5),\; B \mapsto \left( \begin{smallmatrix} B & 0 \\ 0 & B^{-1} \end{smallmatrix} \right) \cdot \Ug(5)$\, is a totally
geodesic embedding of this type.
\end{enumerate}

\subsection[$E_6/(\Ug(1)\cdot\Spin(10))$]{$\boldsymbol{E_6/(\Ug(1)\cdot\Spin(10))}$}
\label{SSe:rk2:EIII}

% See Ref.~\refcite{Klein:2008-tgex}: Section~3.2 for the classification, Section~3.3 for the embeddings.

\begin{enumerate}
\item \,$\OP^2_{\vkap=1/2}$\, 
%This is a real form of the Hermitian symmetric space \,$E_6/(\Ug(1)\cdot\Spin(10))$\,.
\item \,$\CP^5_{\vkap=1} \times \CP^1_{\vkap=1}$\, 
%This is a meridian of \,$E_6/(\Ug(1)\cdot\Spin(10))$\,, corresponding to \,$\SO(10)/\Ug(5)$\,. 
\item \,$G_2^+(\R^{10})_{\mathrm{srr}=\sqrt2}$\, 
%This is a polar and also the corresponding meridian of \,$E_6/(\Ug(1)\cdot\Spin(10))$\,.
\item \,$G_2(\C^6)_{\mathrm{srr}=1}$\, 
\item \,$(G_2(\HH^4)/\Z_2)_{\mathrm{srr}=1}$\, 
%This is a real form of the Hermitian symmetric space \,$E_6/(\Ug(1)\cdot\Spin(10))$\,.
\item \,$\SO(10)/\Ug(5)_{\mathrm{srr}=1}$\, 
%This is a polar of \,$E_6/(\Ug(1)\cdot\Spin(10))$\,.
\end{enumerate}

\subsection[$E_6/F_4$]{$\boldsymbol{E_6/F_4}$}
\label{SSe:rk2:EIV}

% See Ref.~\refcite{Klein:2008-tgex}, Sections~4.2 and 4.3.

\begin{enumerate}
\item \,$\OP^2_{\vkap=1/4}$\,
\item \,$\HP^3_{\vkap=1/4}$\,
\item \,$((\SU(6)/\Sp(3))/\Z_3)_{\mathrm{srr}=1}$\,
\item \,$(\Sph^9_{r=1} \times \Sph^1_{r=\sqrt{3}})/\Z_4$\,
\end{enumerate}

\subsection[$G_2/\SO(4)$]{$\boldsymbol{G_2/\SO(4)}$}
\label{SSe:rk2:G}

% See Ref.~\refcite{Klein:2008-tgex}, Section~5.4.

\begin{enumerate}
\item \,$\SU(3)/\SO(3)_{\mathrm{srr}=\sqrt{3}}$\,
\item \,$(\Sph^2_{r=1} \times \Sph^2_{r=1/\sqrt{3}})/\Z_2$\,
\item \,$\CP^2_{\vkap=3/4}$\,
\item \,$\Sph^2_{r=\tfrac23\,\sqrt{21}}$\,
\end{enumerate}

\subsection[$\SU(3)$]{$\boldsymbol{\SU(3)}$}
\label{SSe:rk2:A2}

% See Ref.~\refcite{Klein:2008-tgex}, Section~4.5.

\begin{enumerate}
\item \,$\SU(3)/\SO(3)_{\mathrm{srr}=1}$\, \\
The Cartan embedding \,$f: \SU(3)/\SO(3) \to \SU(3)$\, is a totally geodesic embedding of this type.
\item \,$(\Sph^3_{r=1} \times \Sph^1_{r=\sqrt{3}})/\Z_2$\,
\item \,$\CP^2_{\vkap=1/4}$\, \\
The Cartan embedding \,$f: \SU(3)/\mathrm{S}(\Ug(2)\times\Ug(1)) \to \SU(3)$\, is a totally geodesic embedding of this type.
\item \,$\RP^3_{\vkap=1/4}$\, 
\end{enumerate}

\subsection[$\Sp(2)$]{$\boldsymbol{\Sp(2)}$}
\label{SSe:rk2:B2}

% See Ref.~\refcite{Klein:2008-tgex}, Section~3.4.

\begin{enumerate}
\item \,$G_2^+(\R^5)_{\mathrm{srr}=1}$\, \\
The Cartan embedding \,$f: \Spin(5)/(\Spin(2)\times\Spin(3)) \to \Spin(5)\cong \Sp(2)$\, is a totally geodesic embedding of this type.
\item \,$\Sp(1) \times \Sp(1)$\, \\
The canonically embedded \,$\Sp(1)\times\Sp(1)\subset\Sp(2)$\, is a totally geodesic submanifold of this type.
\item \,$\HP^1_{\vkap=1/2}$\, \\
The Cartan embedding \,$f: \Sp(2)/(\Sp(1)\times\Sp(1)) \to \Sp(2)$\, is a totally geodesic embedding of this type.
\item \,$\Sph^3_{r=\sqrt{5}}$\, 
\end{enumerate}

\subsection[$G_2$]{$\boldsymbol{G_2}$}
\label{SSe:rk2:G2}

% See Ref.~\refcite{Klein:2008-tgex}, Section~5.2 for the classification, Section~5.3 for the embeddings. 

\begin{enumerate}
\item \,$G_2/\SO(4)_{\mathrm{srr}=1}$\, \\
The Cartan embedding \,$f: G_2/\SO(4) \to G_2$\, is a totally geodesic embedding of this type.
\item \,$(\Sph^3_{r=1} \times \Sph^3_{r=1/\sqrt{3}})/\Z_2$\, 
\item \,$\SU(3)_{\mathrm{srr}=\sqrt{3}}$\, \\
Regard \,$G_2$\, as the  automorphism group of the division algebra of the octonions \,$\OO$\, and fix an imaginary unit octonion \,$i$\,. Then 
the subgroup \,$\Menge{g\in G_2}{g(i)=i}$\, is isomorphic to \,$\SU(3)$\, and a totally geodesic submanifold of this type.
\item \,$\Sph^3_{r=\tfrac23\sqrt{21}}$\,
\end{enumerate}

\section*{Acknowledgments}
This work was supported by a fellowship within the Postdoc-Programme of the German Academic Exchange Service (DAAD).


\begin{thebibliography}{[20]}
\bibitem[1]{Broecker-tom-Dieck:1985}{T.~Br\"ocker, T.~tom Dieck},
{\it Representations of compact Lie groups}\/,
Springer, New York, 1985.

\bibitem[2]{Chen:1987}{B.-Y.~Chen},
{\it A new approach to compact symmetric spaces and applications}\/,
Katholieke Universiteit Leuven, 1987.

\bibitem[3]{Chen/Nagano:totges1-1977}{B.-Y.~Chen, T.~Nagano},
{\it Totally geodesic submanifolds of symmetric spaces,~I}\/,
Duke Math.~J. 44~(1977), 745--755.

\bibitem[4]{Chen/Nagano:totges2-1978}{B.-Y.~Chen, T.~Nagano}, 
{\it Totally geodesic submanifolds of symmetric spaces,~II}\/,
Duke Math.~J.~45 (1978), 405--425.

\bibitem[5]{Ihara:tgcomplex-1967}{S.~Ihara},
{\it Holomorphic imbeddings of symmetric domains}\/,
J.~Math.~Soc.~Japan 19~(1967), 261--302.

\bibitem[6]{Ikawa/Tasaki:2000}{O.~Ikawa, H.~Tasaki},
{\it Totally geodesic submanifolds of maximal rank in symmetric spaces}\/,
Japan.~J.~Math.~26 (2000), 1--29.

\bibitem[7]{Klein:2007-claQ}{S.~Klein},
{\it Totally geodesic submanifolds of the complex quadric}\/,
Differential Geom.~Appl., 26 (2008), 79--96.

\bibitem[8]{Klein:2007-tgG2}{S.~Klein},
{\it Totally geodesic submanifolds of the complex and the quaternionic 2-Grassmannians}\/,
to appear in Trans.~Amer.~Math.~Soc., \texttt{arXiv:0709.2644}.

\bibitem[9]{Klein:2007-Satake}{S.~Klein},
{\it Reconstructing the geometric structure of a Riemannian symmetric space from its Satake diagram}\/,
to appear in Geom.~Dedicata, \texttt{arXiv:0801.4127}.

\bibitem[10]{Klein:2008-tgex}{S.~Klein},
{\it Totally geodesic submanifolds of the exceptional Riemannian symmetric spaces of rank 2}\/,
submitted for publication, \texttt{arXiv:0809.1319}.

\bibitem[11]{Kobayashi:1972}{S.~Kobayashi},
{\it Transformation groups in Differential geometry}\/,
Springer-Verlag Berlin 1972.

\bibitem[12]{Leung:reflective-1974}{D.~S.~P.~Leung},
{\it On the Classification of Reflective Submanifolds of Riemannian Symmetric Spaces}\/,
Indiana Math.~J., 24 (1974), 327--339.

\bibitem[13]{Leung:reflective-errata-1975}{D.~S.~P.~Leung},
{\it Errata: On the Classification of Reflective Submanifolds of Riemannian Symmetric Spaces}\/,
Indiana Math.~J., 24 (1975), 1199.

\bibitem[14]{Leung:reflective-1979}{D.~S.~P.~Leung},
{\it Reflective submanifolds. III. Congruency of isometric reflective submanifolds and corrigenda to the classification of reflective submanifolds}\/,
J.~Diff.~Geom., 14 (1979), 167--177.

\bibitem[15]{Leung:realforms-1979}{D.~S.~P.~Leung},
{\it Reflective submanifolds. IV. Classification of real forms of Hermitian symmetric spaces.}\/,
J.~Diff.~Geom., 14 (1979), 179--185.

\bibitem[16]{Loos:1969-2}{O.~Loos}, % L2
{\it Symmetric spaces II: Compact spaces and classification}\/,
W.~A.~Benjamin Inc., New York, 1969.

\bibitem[17]{Shafarevich:1994}{I.~R.~Shafarevich},
{\it Basic algebraic geometry I}\/,
2nd Edition, Springer, Berlin, 1994.

\bibitem[18]{Wolf:1963-spheres}{J.~A.~Wolf},
{\it Geodesic spheres in Grassmann manifolds}\/,
Illinois J.~Math., 7~(1963), 425--446.

\bibitem[19]{Wolf:1963-elliptic}{J.~A.~Wolf},
{\it Elliptic spaces in Grassmann manifolds}\/,
Illinois J.~Math., 7~(1963), 447--462.

\bibitem[20]{Zhu/Liang:2004}{F.~Zhu, K.~Liang},
{\it Totally geodesic submanifolds of maximal rank in symmetric spaces}\/,
Science in China Ser.~A, 47 (2004), 264--271.

\end{thebibliography}
\end{document}